\documentclass{article}
\usepackage[usenames,dvipsnames]{pstricks}

\usepackage{comment}
\usepackage{pstricks,psfrag}
\usepackage{xcolor, rotating, epic,eepic}
\usepackage{amssymb}
\usepackage{amsmath}
\usepackage{times} 
\usepackage{epsfig,t1enc}
\usepackage{graphicx}
\usepackage{graphics}
\newtheorem{defi}{Definition}
\newtheorem{theorem}{Theorem}
\newtheorem{pro}{Proposition}
\newtheorem{conj}{Conjecture}

\newtheorem{cor}{Corollary}

\def\II{\hbox{1\kern-.2em\hbox{I}}}

\date{}
\title{Criteria for a split real polynomial} 

\author{J.-M. Billiot \and E. Fontenas
\thanks{Laboratoire Jean Kuntzmann, UMR 5224, Universit\'e Grenoble Alpes, 700 avenue centrale,
38041 Domaine Universitaire de Saint-Martin-d'H\`eres, France.
Jean-Michel.Billiot@univ-grenoble-alpes.fr, Eric.Fontenas@univ-grenoble-alpes.fr}}

\begin{document}
\maketitle
\begin{abstract}
In this article, we establish necessary and sufficient conditions for a polynomial of degree $n$ to have exactly $n$ real roots. A complete study of polynomials of degree five is carried out. The results are compared with those obtained using Sturm sequences. 
\end{abstract}
\vspace{1cm}

\par 

\section{Introduction}
The study of polynomials with only real roots is the subject of much research. \cite{Liu07} proposes results concerning sequences of polynomials and obtains sufficient conditions for these polynomials to have only real roots. Their result is based on the notion of interlacing roots. 
Recently, \cite{Gonz18}, \cite{Gonz19} have studied cubic, quartic and quintic polynomials
and have proposed conditions on the coefficients associated with the Sturm sequences which will determine the multiplicities of the real and complex roots as well as the order of the real roots with respect to the multiplicity.

In this article, we give necessary and sufficient conditions on the coefficients of a polynomial of degree $n$ so that this polynomial has $n$ distinct real roots. The idea here is to construct a sequence of polynomials with real roots, introducing the same assumption about the last constant of the polynomial constructed. We take as an example the case of a polynomial of degree three and then explain the general idea. A complete study of fifth-order polynomials is carried out, highlighting the importance of multiple roots to obtain bounds on the coefficients of polynomials.
Our method is then compared with that proposed by Sturm \cite{Stu35}. Sturm's sequences allow us to obtain conditions on the parameters of a polynomial (very effective when the coefficients are known). Unlike the result proposed here, these conditions are not solved for any polynomial. 
\section{Method}
For simplicity's sake (but our method also works in the general case), we are interested in polynomials where the coefficient of the highest monomial is equal to one (otherwise, we simply divide the polynomial by this coefficient).  We also assume that the coefficient of the second-highest monomial is zero: an appropriate change of variable makes this possible. This means that the sum of the roots is zero.
Depending on the context, the notation $P_n$ or $P_n(x)$ will be used to designate the same polynomial of degree $n$ in the variable $x$.\\
We recall the notion of interlaced polynomial roots:

\begin{defi}
Given both polynomials $P$ and $Q$ of order $n$ and $n-1$ respectively and  $\{\alpha_i\}_{1\leq i\leq n}$ and $\{\beta_j\}_{1\leq j\leq n-1}$ be all real roots of $P$ and $Q$ in nonincreasing order respectively. The roots of $Q$ are interlaced with the roots of $P$ if 
$$\alpha_1\leq \beta_1\leq \alpha_2\leq \beta_2\leq \cdots \leq \beta_{n-2}\leq \alpha_{n-1}\leq \beta_{n-1}\leq \alpha_{n}.$$
\end{defi}
\subsection{Discussion for polynomials of degree two and three}
\noindent 1. Consider $P_2(x)=x^2+p/3$. This polynomial has two distinct real roots as soon as $p<0$.\\
\noindent 2. For the polynomial $P_3(x)=x^3+px+q$, by the Euclidean division by $P'_3(x)=3P_2(x)$, it comes then 
$$P_{3}(x)=\frac{x}{3}P^{'}_{3}(x)-R_{1}(x)=xP_{2}(x)- R_{1}(x)$$
with $R_{1}(x)=-(\frac{2}{3}px+q)$. The polynomial $P_3$  has three real roots if and only if $P_2$  has two real roots and if $R_1$ has a real root interlaced with those of  $P_2$. Because of the previous remark, $P_2$  has two real roots if $p<0$.\\
It then remains to fix the constant $q$ so that the root of $R_1$ is interlaced with those of $P_2$. The root of $R_{1}$ is $\beta^{(1)}_{1}=-\frac{3q}{2p}$ and is interlaced with those of $P_{2}$ if $$-\frac{3q}{2p}\in \left]\alpha^{(2)}_{1}=-\left(-\frac{p}{3}\right)^{1/2};\alpha^{(2)}_{2}=\left(-\frac{p}{3}\right)^{1/2}\right[$$ or $$q\in \left]-\frac{2}{3}p\alpha^{(2)}_{1}; -\frac{2}{3}p\alpha^{(2)}_{2}\right[.$$\\
If we call $R_{1}^{0}$ the function defined by  $R_{1}^{0}(x)=R_1(x)+q=-\frac{2}{3}px$,  we have that $$q\in]R_{1}^{0}(\alpha^{(2)}_{1}); R_{1}^{0}(\alpha^{(2)}_{2})[.$$
This means that $R_{1}\left(\alpha^{(2)}_{1}\right)<0$ and that $R_{1}\left(\alpha^{(2)}_{2}\right)>0$. Note here that the result means that the discriminant of $P_3$ namely $\Delta(P_{3})= -(4p^3+27q^2)$ is not negative.
Moreover, if $p<0$, it is well known that only cases of multiple roots for $P_{3}$ correspond to $\beta^{(1)}_{1}= \alpha^{(2)}_{1}$ or $\beta^{(1)}_{1}=\alpha^{(2)}_{2}$ according with the sign of $q$. Remark that this real double root is also a root of $P_{2}$ and $R_{1}$.

\subsection{General idea}
We construct a sequence of polynomials $(P_n)_{n\geq 2}$ so that $\forall n\geq 3,\,P_{n-1}(x)=\frac{1}{n}P^{'}_{n}(x)$:

$$\left\{\begin{array}{l}
 P_2(x)=x^2+c_0/3\\
 P_3(x)= x^3+c_0x+c_1\\
 P_4(x)= x^4+\frac{4!}{ 2!3!}c_0x^2+\frac{4!}{1!3!}c_1x+\frac{4!}{3!0!}c_2\\
P_5(x)=x^5+\frac{5!}{3!3!}c_0x^3+\frac{5!}{3!2!}c_1x^2+\frac{5!}{1!3!}c_2x+\frac{5!}{0!3!}c_3\\
\cdots\\
 P_n(x)=x^n+a_{n-2}x^{n-2}+a_{n-3}x^{n-3}+\cdots+a_0
\end{array}\right.$$
where
$$\forall k\in\{0,\ldots,n-2\},\quad a_k=\frac{n!}{k!3!}c_{n-2-k}.$$

\noindent From a more general point of view, the main result is based on the following result:
\begin{pro}~\\
1. Set $P_n(x)=xQ_{n-1}(x)-R_{n-2}(x)$.
If $Q_{n-1}$ has $n-1$ real roots and if $R_{n-2}$ has $n-2$ real roots interlaced with those of $Q_{n-1}$, then $P_n$ has $n$ real roots.\\
2. For a polynomial of degree $n$, to have $n$ real roots, its derivative must have
$n-1$ distinct real roots.
\end{pro}

Now consider our sequence of polynomials and divide $P_n$ Euclidean by $P_{n-1}(x)=P'_{n}(x)/n$:
$P_n(x)=\frac{x}{n}P'_{n}(x)-R_{n-2}(x)$. $P_n$ has $n$ distinct real roots as soon as its
derivative has $n-1$ distinct real roots and those of $R_{n-2}$ denoted by $\beta_{1}^{(n-2)}\ldots \beta_{n-2}^{(n-2)}$ are interlaced  with those of $P'_{n}$.\\

The coefficient of highest degree of $R_{n-2}$ is $-\frac{2}{n}a_{n-2}$. If $a_{n-2}<0$ (as the sign of the discriminant of $P_{2}$) and if $R_{n-2}$ has $n-2$  real distinct and interlaced roots with $P^{'}_{n}$ then $P_{n}$ has $n$ distinct real roots.
If we write:
$$R_{n-2}(x)=\frac{1}{n-2}(x-\beta^{(1)}_{1})R^{'}_{n-2}(x)-T_{n-4}(x),$$ as $R^{'}_{n-2}(x)=(n-1)R_{n-3}(x)$, $R^{'}_{n-2}$ has $n-3$ distinct real roots. It is therefore necessary to ensure that $T_{n-4}$ has $n-4$ distinct real roots denoted by $\gamma^{(n-4)}_{1}\ldots \gamma^{(n-4)}_{n-4}$ and besides interlaced with those  $\beta_{1}^{(n-3)}\ldots \beta_{n-3}^{(n-3)}$ of $R^{'}_{n-2}$ and the highest degree coefficient
of $T_{n-4}$  is positive and so on.\\

We will specify this construction later when studying polynomials of degree five when we will explain how to choose the last coefficient $a_{0}$ of the polynomial $P_{n}$ iteratively. 
This leads us now to the presentation of our main result.

\section{Main result}
First of all, we present a characteristic property of interlacing roots. Then, we explain how choosing the integration constant $ a_ {0} $ using the extrema of  $P_{n}$ and the remainder of Euclidean division of $P_{n}$ by $P^{'}_{n}$. At last, we propose some necessary and sufficient assumptions on the coefficients to obtain $n$ distinct real roots for $P_{n}$.
\begin{pro}
Let $P_{n}$ be a polynomial and let $$P_{n}(x)=\frac{x}{n}P'_{n}(x)-R_{n-2}(x)=xP_{n-1}(x)-R_{n-2}(x).$$
Denote by $\alpha^{(n-1)}_{1}\ldots\alpha^{(n-1)}_{n-1}$, the $n-1$ distinct real roots of  $P_{n-1}$ and  $\beta^{(n-2)}_{1}\ldots\beta^{(n-2)}_{m}$, the $n-2$ distinct real roots of  $R_{n-2}$.\\
If $n$ is even, $$\alpha_{1}^{(n-1)}<\beta_{1}^{(n-2)}<\alpha_{2}^{(n-1)}<\beta_{2}^{(n-2)}<\ldots <\alpha_{n-2}^{(n-1)}<\beta_{n-2}^{(n-2)}<\alpha_{n-1}^{(n-1)}$$
$$ \Leftrightarrow \sup_{k\in{\{1,\ldots, \frac{n-2}{2}}\}}R_{n-2}(\alpha_{2k}^{(n-1)})<0<\inf_{k\in{\{0,\ldots, \frac{n-2}{2}}\}}R_{n-2}(\alpha_{2k+1}^{(n-1)}).$$\\
If $n$ is odd,
$$\alpha_{1}^{(n-1)}<\beta_{1}^{(n-2)}<\alpha_{2}^{(n-1)}<\beta_{2}^{(n-2)}<\ldots< \alpha_{n-2}^{(n-1)}<\beta_{n-2}^{(n-2)}<\alpha_{n-1}^{(n-1)}$$
$$ \Leftrightarrow \sup_{k\in{\{0,\ldots, \frac{n-3}{2}}\}}R_{n-2}(\alpha_{2k+1}^{(n-1)})<0<\inf_{k\in{\{1,\ldots, \frac{n-1}{2}}\}}R_{n-2}(\alpha_{2k}^{(n-1)}).$$

\end{pro}
\noindent{\bf Proof:} we only present here the proof for $n $
even.
Suppose that $R_{n-2}$ admits $n-2$ distinct roots $\beta_{1}^{(n-2)}\ldots \beta_{n-2}^{(n-2)}$ interlaced with those of $P_{n-1}$ : 
$$\alpha_{1}^{(n-1)}<\beta_{1}^{(n-2)}<\alpha_{2}^{(n-1)}<\beta_{2}^{(n-2)}<\ldots <\alpha_{n-2}^{(n-1)}<\beta_{n-2}^{(n-2)}<\alpha_{n-1}^{(n-1)}.$$
So, 
$$\sup_{k\in{\{0,\ldots, \frac{n-2}{2}}\}}P_{n}\left(\alpha_{2k+1}^{(n-1)}\right)<0<\inf_{k\in{\{1,\ldots, \frac{n-2}{2}}\}}P_{n}\left(\alpha_{2k}^{(n-1)}\right).$$
which is also written
$$\sup_{k\in{\{1,\ldots, \frac{n-2}{2}}\}}R_{n-2}\left(\alpha_{2k}^{(n-1)}\right)<0<\inf_{k\in{\{0,\ldots, \frac{n-2}{2}}\}}R_{n-2}\left(\alpha_{2k+1}^{(n-1)}\right).$$
Conversely, 

$$\begin{array}{c}
\displaystyle\sup_{k\in{\{1,\ldots, \frac{n-2}{2}}\}}R_{n-2}\left(\alpha_{2k}^{(n-1)}\right)<0<\inf_{k\in{\{0,\ldots, \frac{n-2}{2}}\}}R_{n-2}\left(\alpha_{2k+1}^{(n-1)}\right)\\
\Longleftrightarrow \\
\forall k\in\{1,..n-2\},\,R_{n-2}\left(\alpha_{k}^{(n-1)}\right)\times R_{n-2}\left(\alpha_{k+1}^{(n-1)}\right)<0
\end{array}$$
This implies that $R_{n-2}$ has $n-2$ interlacing roots with those of $P_{n-1}$. \\

The following theorem allows to choose $a_0$ so that the new polynomial $P_n$ has $n$ distinct real roots:
\begin{theorem}
Let $P_n$ be a polynomial and set $$P_{n}(x)=\frac{x}{n}P^{'}_{n}(x)-R_{n-2}(x)=xP_{n-1}(x)-R_{n-2}(x).$$
We call $R_{n-2}^0(x)=R_{n-2}(x)+a_0$.\\
Denote by $\alpha^{(n-1)}_{1}\ldots\alpha^{(n-1)}_{n-1}$  the $n-1$ distinct real roots of  $P_{n-1}$.\\

\noindent For $n$ even, if
$$\left\{\begin{array}{l}
\displaystyle\sup_{k\in{\{1,\ldots, \frac{n-2}{2}}\}}R_{n-2}(\alpha_{2k}^{(n-1)})<0<\displaystyle\inf_{k\in{\{0,\ldots, \frac{n-2}{2}}\}}R_{n-2}(\alpha_{2k+1}^{(n-1)})\\
a_{0}\in\left]\displaystyle\sup_{k\in{\{1,\ldots, \frac{n-2}{2}}\}}R_{n-2}^{0}(\alpha_{2k}^{(n-1)});\displaystyle\inf_{k\in{\{0,\ldots, \frac{n-2}{2}}\}}R_{n-2}^{0}(\alpha_{2k+1}^{(n-1)})\right[,
\end{array}\right.$$
then $P_{n}$ has $n$ distinct real roots.\\

\noindent For $n$ odd, if
$$\left\{\begin{array}{l}
\displaystyle\sup_{k\in{\{0,\ldots, \frac{n-3}{2}}\}}R_{n-2}(\alpha_{2k+1}^{(n-1)})<0<\displaystyle\inf_{k\in{\{1,\ldots, \frac{n-1}{2}}\}}R_{n-2}(\alpha_{2k}^{(n-1)})\\
a_0\in\left]\displaystyle\sup_{k\in{\{0,\ldots, \frac{n-3}{2}}\}}R^{0}_{n-2}(\alpha_{2k+1}^{(n-1)})\,;\,\displaystyle\inf_{k\in{\{1,\ldots, \frac{n-1}{2}}\}}R^{0}_{n-2}(\alpha_{2k}^{(n-1)})\right[,
\end{array}\right.$$
then $P_{n}$ has $n$ distinct real roots.
\end{theorem}
The following corollary introduces the notion of degenerate cases in the context of our study. We consider only the case where $n$ is odd (a similar result can be established for $n$ even): 
\begin{cor}\label{cordouble}
In the case of $n$ odd, 
\begin{enumerate}

\item if $a_0=\displaystyle\sup_{k\in{\{0,\ldots, \frac{n-3}{2}}\}}R^{0}_{n-2}(\alpha_{2k+1}^{(n-1)})$ or $a_0=\displaystyle\inf_{k\in{\{1,\ldots, \frac{n-1}{2}}\}}R^{0}_{n-2}(\alpha_{2k}^{(n-1)})$, the polynomial $P_n$ as a double root which is a root of $P_{n-1}$ and $R_{n-2}$.
\item  if $a_{0}=\displaystyle\sup_{k\in{\{0,\ldots, \frac{n-3}{2}}\}}R^{0}_{n-2}(\alpha_{2k+1}^{(n-1)})=\displaystyle\inf_{k\in{\{1,\ldots, \frac{n-1}{2}}\}}R^{0}_{n-2}(\alpha_{2k}^{(n-1)})$, the polynomial $P_n$ has either two double roots or one triple root. So these two real distinct roots are also
roots of $R_{n-2}$ and therefore of $P_n$.
\end{enumerate}
\end{cor}

\noindent Applying recursively Theorem 1, the following theorem allows choosing all coefficients $a_l$, $0\leq l<n-2,$ 
\begin{theorem}\label{teomain}
Let $P_{n}(x)=x^{n}+a_{n-2}x^{n-2}+\ldots+ a_{1}x+a_{0}$ and define the sequence $[P_{i},R_{i},R^{0}_{i}]$, $i\in\{3,\ldots,n\}$, such that 
$$\left\{\begin{array}{l}
P_{i-1}(x)=\displaystyle\frac{1}{i}P'_{i}(x)\\
P_{i}(x)=xP_{i-1}(x)-R_{i-2}(x)\\
\end{array}\right.$$
and, $\forall i\in \{3,\ldots,n\},\, R^{0}_{i-2}(x)=R_{i-2}(x)+a_{n-i}$.\\
We assume that the intervals below of the form $]\sup;\inf[$ have a meaning, i.e., the $\sup$ must always be less than the $\inf$.\\
If $n$ is even, $P_n$ has  $n$ distinct real roots if and only if, for all $l\in\{0,...,\frac{n}2-2\}$,
$$\left\{\begin{array}{l}
\bullet \, a_{n-2}<0\\
\bullet\, a_{2l}\in \left]\displaystyle\sup_{k\in \{1,\ldots, \frac{n-2}{2}-l\}}R^{0}_{n-2l-2}(\alpha^{(n-2l-1)}_{2k})\,;\, \inf_{k\in \{0,\ldots ,\frac{n-2}{2}-1\}}R^{0}_{n-2l-2}(\alpha^{(n-2l-1)}_{2k+1})\right[\\

\bullet \,a_{2l+1}\in\left] \displaystyle\sup_{k\in \{0,\ldots, \frac{n-4}{2}-l\}}R^{0}_{n-2l-3}(\alpha^{(n-2l-2)}_{2k+1})\,;\,  \inf_{k\in \{1,\ldots, \frac{n-2}{2}-1\}}R^{0}_{n-2l-3}(\alpha^{(n-2l-2)}_{2k})\right[.
\end{array}\right.$$
If $n$ is odd, $P_n$ has  $n$  distinct real roots if and only if,\\
$\bullet$ $a_{n-2}<0$.\\
$\bullet$ $\forall\, l\in\{0,...,\frac{n-3}2\},\,
a_{2l}\in\left] \displaystyle\sup_{k\in \{0,\ldots, \frac{n-3}{2}-l\}}R^{0}_{n-2l-2}(\alpha^{(n-2l-1)}_{2k+1})\,;\, \inf_{k\in \{0,\ldots, \frac{n-1}{2}-l\}}R^{0}_{n-2l-2}(\alpha^{(n-2l-1)}_{2k})\right[$\\
$\bullet$ $\forall\, l\in\{0,...,\frac{n-5}2\},
\, a_{2l+1}\in\left]\displaystyle\sup_{k\in \{1,\ldots, \frac{n-3}{2}-l\}}R^{0}_{n-2l-3}(\alpha^{(n-2l-2)}_{2k})\,;\, \inf_{k\in \{0,\ldots, \frac{n-3}{2}-1\}}R^{0}_{n-2l-3}(\alpha^{(n-2l-2)}_{2k+1})\right[.$
\end{theorem}
\noindent{\bf Remark:} The hypothesis that intervals of the form $]\sup;\inf[$ have a meaning will be studied in the context of a polynomial of degree five.

\section{Polynomial of order three and four}
The polynomial $P_3(x)=x^3+px+q$ has three roots if $p<0$ and $q\in \left]-2(-\frac{p}{3})^{3/2}; 2(-\frac{p}{3})^{3/2}\right[\Leftrightarrow\Delta(P_3)=-4p^3-27q^2>0$.

\begin{pro}\label{ProP4} $P_{4}(x)=x^4+2px^2+4qx+4r$ has four distinct real roots if we choose:
$$p<0,\quad q\in \left]-2(-\frac{p}{3})^{3/2}; 2(-\frac{p}{3})^{3/2}\right[,\quad 4r\in ]R_{2}^{0}(\alpha_{2}^{(3)});\inf_{i\in{\{1,3\}}}R^{0}_{2}(\alpha_{i}^{(3)}) [\subset ]R^{0}_{2}(\beta^{(1)}_{1});+\infty[ $$
where 
$$\left\{\begin{array}{l}
\beta^{(1)}_{1}=-\frac{3q}{2p}\,\text{is the root of}\, R_{1}(x)=-\frac23px-q,\\
\alpha_{i}^{(3)},i=1..3,\,\text{the roots of}, P_3(x)=x^3+px+q,\\
R_2^{0}(x)=-px^2-3qx.
\end{array}\right.$$
More precisely, $P_4$ has four distinct roots if 
$$p<0,\quad q\in \left]-2(-\frac{p}{3})^{3/2}; 2(-\frac{p}{3})^{3/2}\right[,\quad 4r\in ]R_{2}^{0}(x_3);R^{0}_{2}(x_1)[$$
where 
$$\left\{\begin{array}{l}
R_{2}^{0}(x_3)=\frac{4p^2}{9}\cos(\theta/3)\left[\cos(\theta/3)+\cos(\theta)\right]\\
R_{2}^{0}(x_1)=\frac{4p^2}{9}\cos(\theta/3+\pi/3)\left[\cos(\theta/3+\pi/3)+\cos(\theta)\right]\\
\cos(\theta)=-\frac{|q|}{2\sqrt{-p/3}}
\end{array}\right.$$
\end{pro}
\noindent{\bf Remarks:}
\begin{enumerate}
\item Say $4r\in ]R^{0}_{2}(\beta^{(1)}_{1});+\infty[$ means $R_2$ has two distinct roots $\beta_{1}^{(2)}$ and $\beta_{2}^{(2)}$ (the discriminant of $R_2$ is given by $\Delta{R_2}=9q^2-16pr>0$). Clearly:
$$\Delta_{2}>0\Leftrightarrow R_{2}(\beta^{(1)}_{1})=R^{0}_{2}(\beta^{(1)}_{1})-4r<0.$$
\item Say $4r\in ]R_{2}^{0}(\alpha_{2}^{(3)});\displaystyle\inf_{i\in{\{1,3\}}}R^{0}_{2}(\alpha_{i}^{(3)}) [$ means the roots of $R_2$ are interlaced with $\alpha_{1}^{(3)}$, $\alpha_{2}^{(3)}$, $\alpha_{3}^{(3)}$, the roots of $P_3$.\\

\end{enumerate}

\noindent{\bf Study of multiple roots for a polynomial of order four}\\
In the proposition \ref{ProP4}, the maximum or minimum values reached by $r$ give very particular cases concerning the roots of the polynomial $P_4$.\\
1. The bound $r=p^2/4$ ($q=0$) corresponds to the fact that the polynomial $P_4$ has two double roots $-\sqrt{-p}$ and $\sqrt{-p}$.\\
2. The bound $r=-p^2/12$ corresponds to a double root $-\sqrt{-p/3}$ for the polynomial $P_3$ (which implies $q^2=-4p^3/27$): so, according with the sign of $q$, $P_4$ has one triple root $\pm\sqrt{-p/3}$.
 
\section{Polynomial of order five}
In this section, we focus on the fifth order. In particular, we explain our assumptions on the choice of the last parameter $s$. Different cases related to multiple roots are studied.

\subsection{Our assumptions}
Let's consider the following polynomial $$P_{5}(x)=x^5+\frac{10p}{3}x^3+10qx^2+20rx+20s=\frac{xP_4(x)}{5}-R_{3}(x)$$ with $R_3(x)=-(\frac43px^3+6qx^2+16rx+20s)$.
For our method, recall the assumptions of order four:
$$p<0,\quad q\in \left]-2\left(-\frac{p}{3}\right)^{3/2};2\left(-\frac{p}{3}\right)^{3/2}\right[,\quad 4r\in \left]R_{2}^{0}(\alpha_{2}^{(3)});\inf_{k\in{\{1,3\}}}R^{0}_{2}(\alpha_{k}^{(3)})\right[ $$
where $\alpha_{i}^{(3)},\,i=1..3,$ are the roots of $P_3(x)=x^3+px+q$ and $R_2^0(x)=-px^2-3qx$.\\
Applying Theorem 1, we explain how to choose the parameter $s$.
 
\subsection{Choice of $s$}
As $P_{5}(x)=\displaystyle\frac{xP_4(x)}{5}-R_{3}(x)$, $R_3(x)$ must have three distinct real roots and
then they have to be interlaced with the four roots $\alpha^{(4)}_{{i}}, i=1\ldots 4,$ of $P_4$. According to the result of order three, we find
$$20s\in ]R^{0}_{3}(\beta^{(2)}_{2}),R^{0}_{3}(\beta^{(2)}_{1})[$$
where $\beta^{(2)}_{i},\, i=1..2$, roots of $R_2(x)=-px^2-3qx-4r$.\\
The Euclidean division of $R_{3}$  by $R^{'}_{3}$ 
give as a remainder whose sign we change
$$T_{1}(x)=\frac{2(16pr-9q^2)}{3p}x+20s-\frac{8qr}{p}$$
and this remainder vanishes at 
$$\gamma^{(1)}_{1}=\frac{3(-10sp+4qr)}{16pr-9q^2}$$
It is enough now that $\gamma^{(1)}_{1}$ is interlaced with
$\beta^{(2)}_{1}$ and $\beta^{(2)}_{2}$:
 $$\beta^{(2)}_{1}< \gamma^{(1)}_{1}< \beta^{(2)}_{2}.$$
We obtain :
$$ 20s\in ]R^{0}_{3}(\beta^{(2)}_{2});R^{0}_{3}(\beta^{(2)}_{1})[=\left] \frac{-9q^3}{p^2} +\frac{24qr}{p} -\frac{ (\Delta_{2})^{3/2}}{3p^2};\frac{-9q^3}{p^2} +\frac{24qr}{p} +\frac{ (\Delta_{2})^{3/2}}{3p^2}\right[.$$ 
With the assumption $\Delta(R_2)=9q^2-16pr>0$,  we deduce that  $R_{3}$ will have three distinct real roots.\\
Or, as $$-T_{1}(\beta^{(2)}_{1})=R_{3}(\beta^{(2)}_{1})= R^{0}_{3}(\beta^{(2)}_{1})-20s$$ and $$-T_{1}(\beta^{(2)}_{2})=R_{3}(\beta^{(2)}_{2})= R^{0}_{3}(\beta^{(2)}_{2})-20s.$$
The condition can be written as $R_{3}(\beta^{(2)}_{1})> 0$ and $R_{3}(\beta^{(2)}_{2})< 0$ or  $T_{1}(\beta^{(2)}_{1})<0$ and $T_{1}(\beta^{(2)}_{2})> 0$. Now, if the three roots of $R_{3}$ are interlaced with those of $P_{4}$, then $$20s\in{ ]\sup_{k\in\{1,3\}}R^{0}_{3}(\alpha^{(4)}_{k});\inf_{k\in\{2,4 \}}R^{0}_{3}(\alpha^{(4)}_{k})[}\subset]R^{0}_{3}(\beta^{(2)}_{2});R^{0}_{3}(\beta^{(2)}_{1})[$$
under the assumption 

$$\sup_{k\in\{1,3\}}R^{0}_{3}(\alpha^{(4)}_{k})<\inf_{k\in\{2,4 \}}R^{0}_{3}(\alpha^{(4)}_{k}).$$

\begin{theorem}\label{teoP5}
Under the condition
\begin{equation}\label{ass}
\sup_{k\in\{1,3\}}R^{0}_{3}(\alpha^{(4)}_{k})<\inf_{k\in\{2,4 \}}R^{0}_{3}(\alpha^{(4)}_{k}).
\end{equation}
 the polynomial $P_5(x)=x^5+\frac{10p}{3}x^3+10qx^2+20rx+20s$ has five distinct real roots if and only if
 $$p<0,\quad q\in \left]-2\left(-\frac{p}{3}\right)^{3/2};2\left(-\frac{p}{3}\right)^{3/2}\right[,\quad 4r\in \left]R_{2}^{0}(\alpha_{2}^{(3)});\inf_{k\in{\{1,3\}}}R^{0}_{2}(\alpha_{k}^{(3)})\right[ $$
 and $$20s\in{ ]\sup_{k\in\{1,3\}}R^{0}_{3}(\alpha^{(4)}_{k});\inf_{k\in\{2,4 \}}R^{0}_{3}(\alpha^{(4)}_{k})[}$$
where 
$$\left\{\begin{array}{l}
\alpha_{i}^{(3)},\,i=1..3,\,\text{are the roots of}, P_3(x)=x^3+px+q\\
R_2^0(x)=-px^2-3qx\\
\alpha_{i}^{(4)},\,i=1..4,\,\text{are the roots of}\, P_4(x)=x^4+2px^2+4qx+4r\\
R_3^0(x)=-\frac43px^3-6qx^2-16rx.\\
\end{array}\right.$$

\end{theorem}

\subsection{Study of the hypothesis of the theorem \ref{teoP5}}

We must then examine the conditions under which the previous hypothesis (\ref{ass}) is satisfied. The assumption (\ref{ass}) may not be satisfied if
$R^{0}_{3}(a)>R^{0}_{3}(b)$ with $a=\alpha^{(4)}_{1}$ the smallest root and  $b=\alpha^{(4)}_{4}$ the biggest root of $P_{4}$ (as the sum of the roots is zero, $a<0$ and $b>0$). Taking $X=ab\leq -2\sqrt{|r|}$, the inequality (\ref{ass}) can be written as  
\begin{equation}\label{ass2}
(\frac{6r}{p}-p)X(X^2-4r)+(\frac{9q^2}{p}+2r)X^2-8r^2<0.
\end{equation}
Under the assumptions
\begin{equation}\label{ass3}
p<0,\quad q\in \left]-2(-\frac{p}{3})^{3/2};2(-\frac{p}{3})^{3/2}\right[,\quad 4r\in ]R_{2}^{0}(\alpha_{2}^{(3)});\inf_{k\in{\{1,3\}}}R^{0}_{2}(\alpha_{k}^{(3)})[,
\end{equation}
$P_4$ has four distinct real roots $a<b<c<d$. Let $Y=a+b$. We have 
$$\left\{\begin{array}{l}
Y+c+d=0\\
(c+d)Y+X+cd=2p\\
(c+d)X+cdY=-4q\\
Xcd=4r\\
\end{array}\right.$$
It follows that 
$$\left\{\begin{array}{l}
-Y^2+X+\frac{4r}{X}=2p\\
-YX+\frac{4r}{X}Y=-4q\\
\end{array}\right.$$
We deduce that $$16q^2X^3=(4r-X^2)^2(X^2-2pX+4r)<0$$ with $r<\frac{p^2}{4}$ and $p-\sqrt{p^2-4r}<X< -2\sqrt{|r|}$.
The condition $R^3(a)<R^3(b)$ is written
$$(X^2-4r)[X(6r-p^2)+2rp]+9q^2X^2>0.$$
Multiplying this inequality by $X$ and noting that
$$16q^2X^3=(4r-X^2)^2(X^2-2pX+4r)<0,$$ we get
\begin{equation}\label{eq1}(X-x_1)(X-x_2)(X-x_3)(X-x_4)< 0
\end{equation}
where
$$\left\{\begin{array}{l}
x_1=\frac{4p+2\sqrt{4p^2-27r}}{3}\\
x_2=\frac{4p-2\sqrt{4p^2-27r}}{3}\\
x_3=\frac{-p+\sqrt{p^2+12r}}{3}\\
x_4=\frac{-p-\sqrt{p^2+12r}}{3}\\
\end{array}\right.$$
From this, we can deduce:
\begin{pro}\label{pro3}
Under the hypothesis (\ref{ass3}), the assumption
\begin{equation}\label{ass6}
\sup_{k\in\{1,3\}}R^{0}_{3}(\alpha^{(4)}_{k})<\inf_{k\in\{2,4 \}}R^{0}_{3}(\alpha^{(4)}_{k})
\end{equation}
is satisfied as soon as:
\begin{enumerate}
\item $r< p^2/9$,
\item $\frac{p^2}{9}\leq r\leq \frac{5p^2}{36}$ and $X< x_1=\frac{4p+2\sqrt{4p^2-27r}}{3}$\\
\item $\frac{5p^2}{36}\leq r\leq \frac{4p^2}{27}$ and $x_2=\frac{4p-2\sqrt{4p^2-27r}}{3}< X<x_1=\frac{4p+2\sqrt{4p^2-27r}}{3}$
 \end{enumerate} 
 with $p<0$ and $$16q^2X^3=(4r-X^2)^2(X^2-2pX+4r)<0.$$
\end{pro}
\noindent{\bf Remarks:} 1. The hypothesis (\ref{ass6}) reduces the interval of values of the parameter $r$, $r\in]-p^2/12;p^2/4[$, obtained by the conditions (\ref{ass3}), conditions for the polynomial $P_4$ to have four real roots.\\
2. The proposition \ref{pro3} also states that for $p<0,\quad q\in \left]-2(-\frac{p}{3})^{3/2};2(-\frac{p}{3})^{3/2}\right[ $, it is not possible to construct a polynomial of degree five with five real roots with $r\in]4p^2/27; p^2/4[$.
\subsection{Notes: degenerate cases}
Three degenerate cases can be distinguished:
\begin{enumerate}
\item $20s=\displaystyle\sup_{k\in\{1,3\}}R^{0}_{3}(\alpha^{(4)}_{k})=R^{0}_{3}(\alpha^{(4)}_{1})=R^{0}_{3}(\alpha^{(4)}_{3})$;
\item $20s=\displaystyle\inf_{k\in\{2,4 \}}R^{0}_{3}(\alpha^{(4)}_{k})=R^{0}_{3}(\alpha^{(4)}_{2})=R^{0}_{3}(\alpha^{(4)}_{4})$;
\item $20s=\displaystyle\sup_{k\in\{1,3\}}R^{0}_{3}(\alpha^{(4)}_{k})=\displaystyle\inf_{k\in\{2,4 \}}R^{0}_{3}(\alpha^{(4)}_{k}).$
\end{enumerate}
These correspond to the fact that the polynomial $P_5$ has two double roots (see corollary \ref{cordouble}). If we denote by $a$, $b$, $c$ the roots of $P_5$ with $a$ and $b$ the double roots, the first two cases correspond to the situations $a<b<c$ or $c<a<b$ and the last case to $a<c<b$. This is the subject of the two following paragraphs.

\subsubsection{Analysis of the first two cases}
For example, here is a polynomial of degree 5 (see Figure \ref{fig1}) corresponding to the first two cases:

\begin{figure}[ht]
\centering
\includegraphics[scale=0.4]{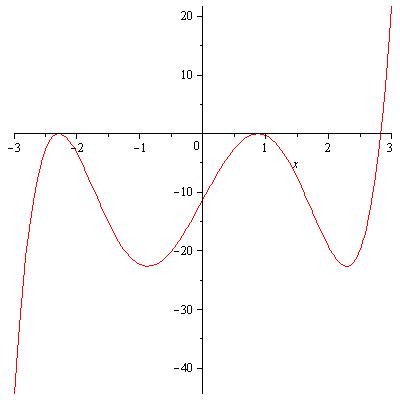}
\caption{an example of degenerate case with $a<b<c$}
\label{fig1}
\end{figure}

~\\
The polynomial $P_5$ has two double roots $a$ and $b$ which are roots of $P_4$ and $R_3$. The Euclidean division of $ P_{5}(x)$ by $ (x-a)^{2} (x-b)^2 $ gives a remainder that is identically zero if
$$\left\{\begin{array}{l}
-3Y^2+2X= \frac{10p}{3}\\
2Y^3+2XY=10q\\
-X(4Y^2-X)=20r\\
2X^2Y=20s
\end{array}\right.
$$
where $X=ab$ and $Y=a+b$.
It comes $$\left\{\begin{array}{l}
Y^2=\frac{2}{3}(X-\frac{5p}{3})\\
Y=\frac{3q}{X-\frac{2p}{3}}\\
X^2-\frac{8p}{3}X+12r=0\\
\end{array}\right.
$$ 
The discriminant of $X^2- \frac{8p}{3}X+12r=0$ is $\Delta=64p^2/9-48r$. The latter is positive for $r<4p^2/27$. Let $x_1= (8p/3+\sqrt{\Delta})/2$ and $x_2= (8p/3-\sqrt{\Delta})/2$ be the two roots of this polynomial.
If $r<0$, $x_1>0$ (reached in r=0) and $x_2<5p/3$. This amounts to $X>0$ either two positive or negative double roots. We can then establish the following result:

\begin{pro}\label{pro5}
The polynomial $P_5$ has two consecutive double roots $a$ and $b$ ($c<a<b$ or $a<b<c$) if
$$\left\{\begin{array}{l}
p<0\\
-\frac{p^2}{12}\leq r\leq \frac{4p^2}{27}\\
x_1\leq -\frac{p}{3}\\
729q^2=-8[11p^3-81pr+2(4p^2-27r)^{3/2}]\\
10s=\frac{3q}{x_1-\frac{2p}{3}}x_1^2
\end{array}\right.$$
The parameter $q$ is in the interval $\left]-\displaystyle\sqrt{\frac{-4p^3}{27}}\,;\,\sqrt{\frac{-4p^3}{27}}\right[$.
\end{pro}

\subsubsection{The last case}
The last case of degeneration corresponds to the following example (see Figure \ref{fig2}):
\begin{figure}[ht]
\centering
\includegraphics[scale=0.3]{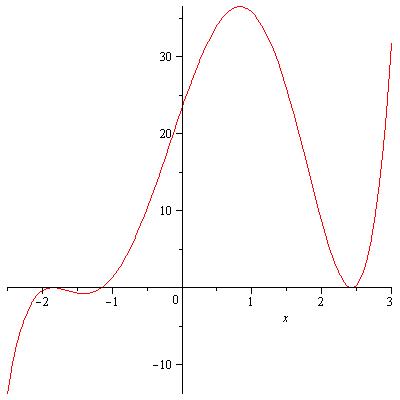}
\caption{an example of the last degenerate case $a<c<b$}
\label{fig2}
\end{figure}

\noindent By using the result of the proposition \ref{pro3}, we obtain the following result
\begin{pro}\label{pro6}
The polynomial $P_5$ has two double roots $a$ and $b$ ($a<c<b$) when
$$\left\{\begin{array}{l}
p<0\\
\frac{5p^2}{36}\leq r\leq \frac{4p^2}{27}\\
x_2=\frac{4p}{3}-\frac{2}{3}\sqrt{4p^2-27r}\\
q^2=-\displaystyle\frac{8}{729}[11p^3-81pr-2(4p^2-27r)^{3/2}]\\
~\\
10s=\frac{3q}{X-\frac{2p}{3}}x_2^2\\
\end{array}\right.$$
The parameter $q$ is in the interval $\left]-\sqrt{\frac{-8p^3}{729}}\,;\,\sqrt{\frac{-8p^3}{729}}\right[$.
\end{pro}
It is possible to find a triplet $(p,q,r)$ for which the propositions \ref{pro5} and \ref{pro6} are satisfied at the same time by choosing $r=4p^2/27$ (i.e. $X=4p/3$), $q^2=-8p^3/729$, $s=4pq/5$.This situation corresponds to a triple root and a double root for $P_5$ of the form $\pm \frac23\sqrt{-2p}$ and $\mp\sqrt{-2p}$.

\section{Comparison with Sturm's method}
Let us first recall Sturm's method. Let $P_n$ be a polynomial and denote its derivative by $P_{n-1}$. Proceed as in Euclid's theorem to find :
$$\left\{\begin{array}{l}
S_{n}=P_{n}\\
S_{n-1}=P^{'}_{n}\\
S_{n}=Q_{1}S_{n-1}-S_{n-2}\\
\ldots\\
S_{2}=Q_{n-1}S_{1}-S_{0}.
\end{array}\right.$$
\begin{theorem} (Sturm's Theorem) The number of distinct real roots of a polynomial in $[a;b]$ is equal to the excess of the number of changes of sign in the sequence $S_n(a),\cdots, S_{0}(a)$ over the number of changes of sign in the sequence $S_n(b),\cdots, S_{0}(b)$.
\end{theorem}
Since we are interested only in polynomials of degree $n$ with $n$ distinct real roots, a corollary of Sturm's theorem is the following result:
\begin{cor}
A polynomial of degree $n$ has $n$ distinct real roots if and only if all the coefficients of the monomial of the greatest degree of the polynomials $S_n(x),\cdots, S_0(x)$ are strictly positive.
\end{cor}
We compare our results with those obtained using Sturm's method for polynomials of degrees 3, 4 and 5.
\subsubsection*{Sturm sequence associated with $P_3$}
In our case, the polynomial $P_3$ has three roots if $p<0$ and $\Delta(P_3)>0$. By Sturm's method,
$$\left\{\begin{array}{l}
S_{3}(x)=P_{3}(x)=x^3+px+q\\
S_{2}(x)=P^{'}_{3}(x)=3x^2+p\\
S_{1}(x)=R_{1}(x)=-(\frac{2}{3}px+q)\\
S_{0}=\frac{-4p^3-27q^2}{4p^2}=\frac{\Delta(P_3)}{4p^2}
\end{array}\right.$$
\subsubsection*{Sturm sequence associated with $P_4$}
Sturm's sequence is given by
$$\left\{\begin{array}{l}
S_{4}(x)=P_{4}(x)=x^4+2px^2+4qx+4r\\
S_{3}(x)=P^{'}_{4}(x)=4x^3+4px+4q\\
S_{2}(x)=R_{2}(x)=-px^2-3qx-4r\\
S_{1}(x)=\displaystyle\frac{-4}{p^2}[(-4pr+p^3+9q^2)x +q(12r+p^2)]\\
S_{0}=\displaystyle\frac{p^2\Delta(P_{4})}{256(-4pr+p^3+9q^2)^2}
\end{array}\right.$$
Using Sturm's theorem, we have four distinct real roots if the coefficients of the term of the highest
degree of $S_{3}, S_2, S_1$ and $S_0$ are strictly positive. These assumptions are the following: $$p<0,\,-4pr+p^3+9q^2<0,\,\Delta(P_{4})>0.$$
In contrast to the results of Sturm's algorithm, we obtained the exact solutions for a polynomial of degree 4 to have exactly four distinct real roots (see Proposition \ref{ProP4}).
\subsubsection*{Sturm sequence associated with $P_5$}
Sturm's sequence is given by
$$\left\{\begin{array}{l}
S_5(x)=P_{5}(x)=x^5+\frac{10p}{3}x^3+10qx^2+20rx+20s\\
S_{4}(x)=P^{'}_{5}(x)=5x^4+10px^2+20qx+20r\\
S_{3}(x)=R_{3}(x)=-(\frac{4}{3}px^3+6qx^2+16rx+20s)\\
S_{2}(x)=-5\left[\frac{(8p^3+81q^2-48pr)x^2}{4p^2}+\frac{(-15ps+4p^2q+54qr)x}{p^2}+4r+\frac{135qs}{2p^2}\right]\\
S_1(x)=-\frac{32p^2}{3(8p^3-48pr+81q^2)^2}(a_{S_1}\,x +b_{S_1})\\
S_{0}=\frac{1000000}{81}\left(\frac{8p^3-48pr+81q^2}{p\,a_{S_1}}\right)^2\Delta(P_{5})
\end{array}\right.$$
with \\
$\left\{\begin{array}{lll}
a_{S_1}&=&-80p^4r-2106q^2pr+1056p^2r^2-3456r^3+240p^2qs+3240qsr\\
&&+40p^3q^2+729q^4-450ps^2\\
b_{S_1}&=&-120p^4s-1755spq^2+1560p^2rs-4320r^2s+40p^3qr+729q^3r\\
&&-864qpr^2+2025s^2q.
\end{array}\right.$\\

\noindent The result of Sturm gives five distinct real roots for $P_{5}$ under the following assumptions:
$$p<0,\,8p^3+81q^2-48pr<0,\,a_{S_1}<0,\Delta(P_5)>0.$$

\subsubsection*{Conjecture}
By observing the last term $S_0$ of the preceding Sturm sequences, we can conjecture the following result

\begin{conj}
Let $P_n$ be the polynomial of order $n$, $P_{n-1}=P'_n/n$, such that :
$$\left\{\begin{array}{l}
P_{n}(x)=xP_{n-1}(x)-R_{n-2}(x)\\
R_{n-2}(x)=R^{0}_{n-2}(x)-a_{0}.
\end{array}\right.$$
The last term $S_0$ given in the Sturm sequence associated with $P_n$ is equal to:
$$
\begin{array}{lll}
S_{0}&=&K_n\Delta(P_{n})=K(-1)^{n(n-1)/2}n^n\displaystyle\prod_{i=1}^{n-1}P_{n}(\alpha^{(n-1)}_{i})
\\
&=&K_n(-1)^{(n-1)(n+2)/2}n^n\displaystyle\prod_{i=1}^{n-1}[R^{0}_{n-2}(\alpha_{i}^{(n-1)})-a_0]\end{array}$$
 where $\alpha^{(n-1)}_{i}, i=1,\ldots, n-1$, distinct real roots of $P_{n-1}$ and $K_n$ strictly not negative constant. $\Delta(P_{n})$ represents the discriminant of $P_{n}$.
\end{conj}
The conjecture concerns the first equality $S_{0}=K_n\Delta(P_{n})$ where $K$ is a strictly positive constant. The other equalities arise from the relationships between the discriminant of a polynomial and the resultant associated with $P_n$ and its derivative $P'_{n}$ (see for example \cite{Gelf94} for results on resultants and Sylvesters' matrices). As we are looking for conditions for $P_n$ to have $n$ distinct real roots, we know that $pgcd(P_n, P'_n)$ is a constant. However, it has not been established, using Euclid's algorithm for example (or by constructing the Sturm sequence associated with $P_n$ in a similar way), that the latter is equal to $K\Delta(P_{n})$.

\section{Concluding remarks}

Our assumptions are explicit and depend on the roots of the previous order. That is why after order five, things become harder. Indeed, the exact expressions of the roots are unknown. Obviously, for cases of orders three and four, Cardano's, Ferrari's,  Descartes's or Euler's formula of the roots are available.\\ 
Otherwise,  Cayley  \cite{Cayley61} and more recently  \cite{Dum91} and \cite{Laval05} proposed different methods for some kind of polynomials of degree six of those roots as functions of the roots for solvable quintics.
For polynomials of order six, we may ask how the assumption we take in theorem 1 can be written using two double roots in particular with the product of the first and the fourth roots or the second and the fifth roots of the polynomial of order five.

Moreover, the conjecture we present for the last term in Sturm's sequence raises many questions. For example, what is the link between discriminants of polynomials of different degrees?
Besides notice the important role played by multiple roots and possible use to obtain upper and lower bounds for intervals of parameters for polynomials having only real roots.  

 On the other side, Ramanujan \cite{Bern94} solve some polynomials of degree three, four, five, six and seven. His approach is very original. It seems that he often started with roots having product one. In this paper, we assume that the sum of the roots is zero even if our result is true in the general case, but it takes a more complex form.

\end{document}